\documentclass[12pt]{article}

\usepackage{amsmath}
\usepackage{amssymb}
\usepackage{latexsym}

\newcommand{\mn}{\medskip\noindent}
\newcommand{\bn}{\bigskip\noindent}
\newcommand{\sn}{\smallskip\noindent}

\newcommand{\cD}{{\mathcal{D}}} 
\newcommand{\Hh}{{\mathcal{H}}}
\newcommand{\cB}{{\mathcal{B}}}
\newcommand{\cN}{{\mathcal{N}}}
\newcommand{\cC}{{\mathcal{C}}}
\newcommand{\cA}{{\mathcal{A}}}
\newcommand{\cU}{{\mathcal{U}}}
\newcommand{\cY}{{\mathcal{Y}}}

\newcommand{\cZ}{{\mathcal{Z}}}

\newcommand{\dN}{\mathbb{N}}
\newcommand{\dC}{\mathbb{C}}
\newcommand{\dR}{\mathbb{R}}
\newcommand{\dB}{\mathbb{B}}

\begin{document}
\vspace{10.6mm}
\title{}\date{}
\begin{center}
{\bf {\Large{On the moment problem of closed semi-algebraic sets}}}

\bn
{\large Konrad Schm\"udgen}

\sn
\small{Fakult\"at f\"ur Mathematik und Informatik\\ Universit\"at Leipzig, 
Augustusplatz 10, 04109 Leipzig, Germany\\ E-mail: schmuedg@mathematik.uni-leipzig.de }
\end{center}
\vspace*{1ex}


\mn
{\bf 1. Introduction and main result} 

\sn
The moment problem for compact semi-algebraic sets has been solved in [S1] (see [PD], Section 6.4, for a refinement of the result). In the terminology explained below, this means that each defining sequence $f$ of a compact semi-algebraic set $K_f$ has property $(SMP)$. On the other hand, for many non-compact semi-algebraic sets (for instance, sets containing a cone of dimension two [KM], [PS]) the moment problem is not solvable. Only very few non-compact semi-algebraic sets (classes of real algebraic  curves [So], [KM], [PS] and cylinder sets [Mc]) are known to have a positive solution of the moment problem. 

In this paper we study semi-algebraic sets $K_f$ such that there exist polynomials $h_1,{\dots},h_n$ which are {\it bounded} on the set $K_f$. Our main result (Theorem 1) reduces the moment problem for the set $K_f$ to the moment problem for the ``fiber sets" $K_f\cap C_\lambda$, where $C_\lambda$ is real algebraic variety $\{x\in\dR^d:h_1(x)=\lambda_1,{\dots},h_n(x)=\lambda_n\}$. From this theorem new classes of {\it non-compact} closed semi-algebraic sets are obtained for which the moment problem has an affirmative solution. Combined with a result of V. Powers and C. Scheiderer [PS], it follows that tube sets around certain real algebraic curves have property $(SMP)$ (see Theorem 9).

Let $f=(f_1,{\dots},f_k)$ be a finite set of polynomials $f_j\in\dR[x]\equiv\dR
[x_1,{\dots},x_d]$ and let $K\equiv K_f$ be the associated closed semi-algebraic subset 
defined by
$$
K_f=\{x\in \dR^d:f_1(x)\ge 0,{\dots},f_k(x)\ge 0\}.
$$
Let $T_f$ denote the corresponding preorder, that is, $T_f$ is the set of all finite sums of elements $f_1^{\varepsilon_1}\cdots 
f_k^{\varepsilon_k} g^2$, where $g\in\dR[x]$ and $\varepsilon_1,{\dots},\varepsilon_k
\in\{0,1\}$. We abbreviate $\dC[x]=\dC[x_1,{\dots},x_d]$.

In this paper we investigate the following properties $(MP)$ and $(SMP)$ of the sequence $f$:\\
$(MP)$: {\it For each linear functional $L$ on $\dC [x]$ such that $L(T_f)\ge 0$ 
there is \\
${}$\hspace{1.1cm} a positive Borel measure $\mu$ on $\dR^d$ such that }
\begin{equation}\label{1}
L(p)=\int p(\lambda) d\mu(\lambda)~ for~ all ~ p\in\dC[x]~.
\end{equation}
$(SMP)$:{\it For each linear functional $L$ on $\dC[x]$ such that $L(T_f)\ge 0$ 
there exists \\
${}$\hspace{1.3cm} a positive Borel measure $\mu$ on $\dR^d$ such that ${\rm supp}~\mu\subseteq K_f$ 
and (1)\\
 ${}$\hspace{1.3cm} holds.}\\
By Theorem 1 in [S1], each defining sequence $f$ of a compact semi-algebraic set $K_f$ has property $(SMP)$. However, for non-compact semi-algebraic sets $K_f$ property $(SMP)$ depends in general on the defining sequence (see Example 5 below). Nevertheless, there are also classes of non-compact semi-algebraic sets $K_f$ (see Corollary 9 and Example 4) such that property $(SMP)$ holds for each defining sequence $f$.

We now state the main result of this paper. Suppose that $h_1,{\dots}, h_n\in\dR[x]$ are polynomials which are {\it bounded} on the set $K_f$. Let 
\begin{align*}
&m_j:= \inf \{h_j(x); x\in K_f\}, ~M_j:=\sup \{h_j(x);x\in K_f\},\\
&\Lambda :=[m_1, M_1]\times{\cdots}\times [m_n, M_n]\subseteq \dR^n.
\end{align*}
For $\lambda=(\lambda_1,{\dots},\lambda_n)\in \Lambda$, let $f(\lambda)$ be the sequence of polynomials given 
by
$$
f(\lambda) := (f_1,{\dots}, f_k,h_1{-}\lambda_1,{-}(h_1{-}\lambda_1),{\dots}, h_n{-}\lambda_n,{-}(h_n{-}\lambda_n)).
$$
Clearly, the corresponding semi-algebraic set $K_{f(\lambda)}$ is the intersection of $K_f$ with the algebraic variety 
$$
C_\lambda=\{ x\in\dR^d:h_1(x)=\lambda_1,{\dots}, h_n(x)=\lambda_n\}~.
$$

\mn
{\bf Theorem 1.} {\it Retain the preceding assumptions and notations. Suppose that for all $\lambda\in \Lambda$ the sequence $f(\lambda)$ has property $(MP)$ (resp. $(SMP)$. Then the sequence $f$ has property $(MP)$ (resp. $(SMP)$).}

\mn
The proof of Theorem 1 will be completed at the end of Section 3. The technical ingredients of the proof might be of interest in themselves. Let $L$ be a linear functional on the algebra $\dC[x]$ such that $L(T_f)\ge 0$ and let $\pi_L$ be the corresponding $GNS$-representation. In Section 2 we show that the operator $\pi_L(p)$ is bounded if the polynomial $p\in\dC[x]$ is bounded on the set $K_f$. In Section 3 we decompose the closure of the representation $\pi_L$ into a direct integral of representations $\pi_\lambda,\lambda\in\Lambda$, such that $\pi_\lambda (h_j)=\lambda_j{\cdot} I$ for $j=1,{\dots},n$ almost everywhere. Section 4 contains some applications and examples. As a by-product, we construct a simple explicit example of a positive linear functional on $\dC[x_1,x_2]$ which is not a moment functional.

Let us recall the definition of a $\ast$-representation (see [S2] for more details). Let $\cA$ be a complex unital $\ast$-algebra and let $\cD$ be a pre-Hilbert space with scalar product $\langle\cdot,\cdot\rangle$. A $\ast$-{\it representation} of $\cA$ on $\cD$ is an algebra homomorphism $\pi$ of $\cA$ into the algebra $L(\cD)$ of linear operators mapping $\cD$ into itself such that $\pi(1)=I$ and $\langle\pi(a)\varphi,\psi\rangle=\langle\varphi,\pi(a^\ast)\psi\rangle$ for all $a\in\cA$ and $\varphi,\psi\in\cD$. Here 1 is the unit element of $\cA$ and $I$ denotes the identity map of $\cD$. The closure of a closable Hilbert space operator $A$ is denoted by $\bar{A}$. If $T$ is a subset of $\dC[x]$ and $L$ is a functional on $\dC[x]$, we write $L(T)\ge 0$ if $L(t)\ge 0$ for all $t\in T$. As usual, $\dB(\Hh)$ is the $C^\ast$-algebra of bounded linear operators on the Hilbert space $\Hh$.

\mn
Acknowledgement. I would like to thank C. Scheiderer for  valuable and inspiring correspondence concerning his results.

\mn
{\bf 2. Bounded polynomials on semi-algebraic sets}

\sn
We denote by $\cB(K_f)$ the set of all polynomials in $\dC[x]$ which are bounded on the semi-algebraic set $K_f$. Clearly, $\cB(K_f)$ is a $\ast$-algebra with involution $p(x)\rightarrow\overline{p(x)}.$

For $p\in\cB(K_f)$ we define $\| p\|_K=\sup \{|p(\lambda)|;\lambda\in K_f\}$. Obviously $\cN_f:=\{ p:\|p\|_K=0\}$ is a $\ast$-ideal of $\cB(K_f)$ and $\| p+\cN_f\|:=\| p \|_K, p\in\cB(K_f)$, defines a norm on the quotient $\ast$-algebra $\cB(K_f)/\cN_f$. Let $\cC(K_f)$ denote the completion of $\cB(K_f)/\cN_f$ with respect to this norm. Then $\cC(K_f)$ is an abelian $C^\ast$-algebra which contains $\cB(K_f)/\cN_f$ as a dense $\ast$-subalgebra. If the set $K_f$ is compact, then $\cB(K_f)=\dC[x]$ and $\cC(K_f)$ is the $C^\ast$-algebra of all continuous functions on $K_f$. At the other extreme, if $K_f=\dR^d$, then $\cB(K_p)=\dC {\cdot}I$ and so 
$\cC(K_f)=\dC {\cdot} I$. 

{\it Throughout this paper we suppose that $L$ is a linear functional on\break 
$\dC[x_1,{\dots},x_d]$ such that $L(T_f)\ge 0$.}

Since $p\bar{p}\in T_f$ and hence $L(p\bar{p})\ge 0$ for $p\in\dC[x]$, $L$ is a positive linear functional on the $\ast$-algebra $\dC[x]$ and the GNS-construction applies. To fix the notation, we briefly repeat this construction. The set $\cN=\{p\in\dC[x]:L(p\bar{p})=0\}$ is an ideal of the algebra $\dC[x]$. Hence there is a scalar product $\langle \cdot,\cdot\rangle$ on the quotient space $\cD_L:=\dC[x]/\cN$ and a linear mapping $\pi_L(p):\cD_L\rightarrow \cD_L$ for $p\in\dC[x]$ defined by $\langle p_1{+}\cN, p_2{+}\cN\rangle=L(p_1\overline{p_2})$ and $\pi_L(p)(q+\cN)=(pq+\cN)$ for $p, p_1,p_2,q\in \dC[x]$. For notational simplicity we shall write $p$ instead of $p+\cN$. The map $p{\rightarrow}\pi_L(p)$ is a $\ast$-representation of the $\ast$-algebra $\dC[x]$ on the domain $\cD_L$ such that
\begin{equation}\label{2}
\langle \pi_L(p) q_1,q_2\rangle=L(pq_1\overline{q_2}),~p,q_1,q_2\in\dC[x].
\end{equation}
Let $\Hh_L$ denote the Hilbert space completion of the pre-Hilbert space $\cD_L$.

The next proposition is crucial in what follows. 

\mn
{\bf Proposition 2.} {\it If $p\in\cB(K_f)$, then the operator $\pi_L(p)$ is bounded on $\cD_L$ and $\|\pi_L(p)\|\le\| p\|_K$.}

\mn
{\bf Proof.} The proof essentially repeats arguments used in the proof of Theorem 1 in [S1]. First we note that it suffices to prove the assertion for real polynomial $p\in\cB(K_f)$, because $\| \pi_L (q) \|^2 {=} \| \pi_L (q \bar{q}) \| {\le} \| q \bar{q}\|_K {=} \| q \|^2_K$. Fix $\varepsilon {>} 0$ and put $\rho {=} \| p \|_K {+} \varepsilon$. Then the polynomial $\rho^2 {-} p^2$ is positive on the set $K_f$. Therefore, by G. Stengle's Positivstellensatz ([Sg], see also [PD]), there exist polynomials $h,g\in T_f$ such that $(\rho^2 {-} p^2)g {=} 1 {+} h$.

Let $q\in\dC[x]$. Since the preorder $T_f$ is closed under multiplication, \break
$p^{2n-2}(1+h)q\bar{q}\in T_f$ and hence $L(p^{2n-2}(1+h)q\bar{q})\ge 0$ for $n\in\dN$. Therefore,
\begin{eqnarray}
\quad L(p^{2n}~g q \bar{q})&=& L(p^{2n-2}(\rho^2 g{-}1{-}h)q \bar{q}) \hspace{2.9cm} \notag\\
&=&\rho^2 L(p^{2n-2} g~q \bar{q}) {-} L(p^{2n-2}(1{+}h)q \bar{q}) {\le} \rho^2 L(p^{2n-2} g~q \bar{q})\notag
\end{eqnarray}
and so
\begin{align}\label{3}
&&\quad L(p^{2n} g~q\bar{q})\le \rho^{2n} L(g~q\bar{q}), n\in\dN.\hspace{2.3cm}
\end{align}
As in [S1], we now use the Hamburger moment problem. Define a linear functional $L_1$ on the polynomial algebra $\dC[t]$ by 
$L_1(t^n)=L(p^n q\bar{q}), n\in\dN_0$. Since $L_1(r(t)\bar{r}(t))=L(r(p)\bar{r}(p) q\bar{q})\ge 0$ for $r\in\dC[t]$, there exists a positive Borel measure $\nu$ on the real line such that $L_1(t^n)=\int t^nd\nu(t),n\in \dN_0$. For $\lambda>0$, let $\chi_\lambda$ be the characteric function of the set $\{t\in\dR :|t|>\lambda\}$. Using the fact that $L(p^{2n}(h+p^2 g)q\bar{q})\ge 0$ and inequality (\ref{3}) 
we get 
\begin{align*}
&\lambda^{2n} \int\chi_\lambda d\nu\le \int t^{2n} d\nu (t)=L_1(t^{2n})=L(p^{2n} q\bar{q})\\
&\le L(p^{2n}(1+h+p^2 g)q\bar{q})=L(p^{2n}(\rho^2 g)q\bar{q})\le \rho^{2n+2} L(g~q\bar{q})
\end{align*}
for all $n\in\dN$. If $\lambda >\rho$, the latter implies that $\int \chi_\lambda d\nu=0$, so we have ${\rm supp} \nu\subseteq [-\rho,\rho]$. Hence, by (\ref{2}), 
$$
\|\pi_L(p)q\|^2=L(p^2q\bar{q})=\int t^2 d\nu\le \rho^2 \int d\nu =\rho^2 L(q \bar{q})=\rho^2\| q\|^2,
$$
that is, $\pi_L(p)$ is bounded and $\|\pi_L (p)\|\le \rho=\| p\|_K+\varepsilon$. Letting $\varepsilon\rightarrow 0$, the assertion follows.\hfill$\Box$

\bn
{\bf Corollary 3.} {\it Let $a,b\in\dR$ and $p\in\dR[x]$. If $a\le p\le b$ on the set $K_f$, then $a{\cdot} I\le \pi_L(p)\le b{\cdot} I$ on the Hilbert space $\Hh_L$.}

\mn
{\bf Proof.} Putting $q=p-\frac{a+b}{2}$, we have $\| q\|_K\le \frac{b-a}{2}$ and hence $\|\pi_L (p)\|\le \frac{b-a}{2}$ by Proposition 2 which implies the assertion. \hfill $\Box$

\mn
Each operator $\pi_L(p),p\in \cB (K_f)$, extends by continuity to a bounded linear operator, denoted again by $\pi_L(p)$, on the Hilbert space $\Hh_L$. Since $\|\pi_L (p)\|\le \| p\|_K$, the map $\cB(K_f)\ni p\rightarrow\pi_L(p)\in\dB(\Hh_L)$ passes to a continuous $\ast$-homomorphism of the quotient $\ast$-algebra $\cB(K_f)/\cN_f$ into $\dB(\Hh_L)$. Hence it extends by continuity to a $\ast$-homomorphism of the $C^\ast$-algebra $\cC(K_f)$ into $\dB(\Hh_L)$.

\mn
{\bf Corollary 4.} Let $p\in\dC[x]$. If $p=0$ on $K_f$, then $L(p)=0$ and $L(p \bar{p})=0$, that is, $p\in\cN$.

\mn
{\bf Proof.} By Proposition 2, $\pi_L(p)=0$. Hence, by (2), $L(p)=\langle \pi_L(p)1,1\rangle =0$ and $L(p \bar{p})=\langle \pi_L(p)1,p\rangle=0$.\hfill $\Box$

\sn
Corollary 4 follows also directly from the real Nullstellensatz combined with the Cauchy-Schwarz inequality.

\bn
{\bf 3. A direct integral decomposition of the GNS representation}

\sn
In this section we need some more notions on unbounded representations and on direct integrals. For the former we refer to the monograph [S2]. Especially, we use the closure $\overline{\pi_L}$ of a $\ast$-representation $\pi_L$ (see [S2], 8.1) and direct integrals of measurable family $\lambda{\rightarrow}\pi_\lambda$ of $\ast$-representations ([S2], 12.3). Our main reference on direct integrals is [D], Chapter II.

Let $h_j,m_j,M_j,j{=}1,{\dots},n$, and $K_f,\Lambda$ be as in Section 1. Since the polynomial $h_j$ is bounded on $K_f$, the operator $\pi_L(h_j)$ is bounded by Proposition 2. Let $H_j{:=}\overline{\pi_L(h_j)}$. Then $H{:=}(H_1,{\dots}, H_n)$ is an $n$-tuple of commuting bounded self-adjoint operators on the Hilbert space $\Hh_L$. Let $\cC$ denote the abelian unital $C^\ast$-algebra generated by the operators $H_1,{\dots}H_n$ and $I$. Because a character $\chi$ of $\cC$ is determined by its values on $H_1,{\dots}, H_n$, the map $\chi{\rightarrow}(\chi(H_1),{\dots},\chi(H_n))$ is a homomorphism of the Gelfand spectrum of $\cC$ on a compact subset, denoted by $\sigma(H)$, of $\dR^n$. Since $m_j{\cdot} I\le H_j\le M_j\cdot I$ by Corollary 3, we have $\sigma(H)\subseteq\Lambda =[m_1,M_1]\times {\cdots}\times[m_n,M_n]$. The bicommutant $\cC^{\prime\prime}$ of $\cC$ is the von Neumann algebra generated by $H_1,{\dots},H_n,I$. The next proposition gives a direct integral decomposition of the closure $\overline{\pi_L}$ of the $\ast$-representation $\pi_L$ of the $\ast$-algebra $\dC[x_1,{\dots},x_d]$.

\mn
{\bf Proposition 5.} {\it There exist a regular positive Borel measure $\nu$ on the compact set $\sigma(H)$, a measurable field $\lambda{\rightarrow}\Hh_\lambda$ of non-zero Hilbert spaces $\Hh_\lambda$, a measurable field $\lambda{\rightarrow}\pi_\lambda$ of closed $\ast$-representations of $\dC[x_1,{\dots},x_d]$ and an isometry $U$ of $\Hh_L$ on the Hilbert space $\Hh=\int^\oplus\Hh_\lambda d\nu(\lambda)$ such that:\\
(i) $U\cC^{\prime\prime} U^{-1}$ is the algebra of bounded diagonalizable operators on $\Hh$.\\
(ii) $U\overline{\pi_L} U^{-1}=\int^\oplus \pi_\lambda d\nu (\lambda)$.\\
(iii) $\pi_\lambda (h_j)=\lambda_j{\cdot} I,j=1,{\dots},n, ~\nu{-}a.e.$ on $\sigma(H)$. }

\mn
{\bf Proof.} First we note that Hilbert space $\Hh_L$ is separable, because the linear span of vectors $\pi_L (x^n)1, n\in\dN^d_0$, is dense in $\Hh_L$. Now we apply Theorems 1 and 2 on pages 217 and 220, respectively, in [D] to the $C^\ast$-algebra 
$\cY{=}\cC$ and the von Neumann algebra $\cZ{=}\cC^{\prime\prime}$. We identify the Gelfand spectrum of $\cC$ with the compact set $\sigma(H)$. Then, by these theorems, there is a positive regular Borel measure $\nu$ on $\sigma(H)$, a measurable family $\lambda{\rightarrow}\Hh_\lambda$ of non-zero Hilbert spaces and an isometry $\cU$ of $\Hh_L$ on $\Hh{=}\int^\oplus \Hh_\lambda d\nu(\lambda)$ such that (i) holds. Then, by construction, $H_j{=}\overline{\pi_L(h_j)}{=}\int^\oplus\lambda_j{\cdot} I d\nu(\lambda)$ which yields (iii).

Let $j\in\{1,{\dots}, n\}$. For $p\in\dC[x]$ and $\varphi\in\cD_L$, we have $H_j\pi_L(p)\varphi=\pi_L(h_j)\pi_L(p)\varphi=\pi_L(p)\pi_L(h_j)\varphi=\pi_L(p)H_j\varphi$ and hence $H_j\overline{\pi_L(p)}\subseteq\overline{\pi_L(p)} H_j$. The latter means that $H_j$ belongs to the symmetrized strong commutant $\pi_L(\dC[x])^\prime_{ss}$ (see [S2], Definition 7.2.7). Since $\pi_L(\dC[x])^\prime_{ss}$ is a von Neumann algebra, $\cC^{\prime\prime}\subseteq\pi_L(\dC[x])^\prime_{ss}$. Therefore, we can continue as in the proofs of Theorems 12.3.1 and 12.3.5 in [S2] (note that condition (HS) theorem is satisfied by remark 2 on p. 343) and obtain a direct integral decomposition as in (ii). \hfill$\Box$

\medskip
For notational simplicity we identify $\Hh_L$ and $\Hh$ via the isometry $\cU$. Let $\cD_\lambda$ be the domain of the representation $\pi_\lambda$ and let $\lambda{\rightarrow}1_\lambda$ be the vector field of $\Hh$ corresponding to the unit element $1\in\cD_L$. By Proposition 5, there is a $\nu$-null set $N\subseteq\sigma(H)$ such that $1_\lambda\in\cD_\lambda$ and $\pi_\lambda(h_j)=\lambda_j{\cdot} I, j=1,{\dots},n$, for all $\lambda\in\sigma(h)\backslash N$. For $\lambda\in\sigma(H)/N$ we define a linear functional $L_\lambda$ on $\dC[x]$ by $L_\lambda(p)=\langle\pi_\lambda(p)1_\lambda,1_\lambda\rangle, p\in\dC[x]$. For $\lambda\in N$ set $L_\lambda\equiv 0$.

\mn
{\bf Lemma 6.} (i) $L(q(h)p)=\int_{\sigma(H)} q (\lambda)L_\lambda(p)d\nu(\lambda)$ for $p,q\in\dC[x]$.\\
(ii) {\it There is a $\nu$-null set $N^\prime\supseteq N$ of $\sigma(H)$ such that $L_\lambda(T_f)\ge 0$ and\\ 
$L_\lambda((h_j-\lambda_j)p)=0$ for $\lambda\in\sigma(H)/N^\prime$, $p\in\dC[x], j=1,{\dots},n$.}

\mn
{\bf Proof.} (i): Since $\overline{\pi_L}{=}\int^\oplus \pi_\lambda d\nu(\lambda)$, from (2) that 
\begin{align*}
&L(q(h)p){=}\langle\overline{\pi_L} (q(h)p)1,1\rangle{=}\int\langle \pi_\lambda(q(h)p)1_\lambda, 1_\lambda\rangle d\nu(\lambda){=}\\
&\int \langle \pi_\lambda (q(\lambda)p)1_\lambda, 1_\lambda\rangle d\nu(\lambda){=}\int_{\sigma(H)} q(\lambda) L_\lambda (p) d\nu (\lambda)~.
\end{align*}
(ii): Let $p\in T_f$. For $q\in\dC[x]$ we have $\bar{q}(h)q(h)p\in T_f$ and hence 
$$
L(\bar{q}(h)q(h)p){=}\int_{\sigma(H)} |q(\lambda)|^2 L_\lambda(p) d\nu (\lambda)\ge 0
$$ 
by (i). From the Weierstrass approximation theorem we therefore conclude that $L_\lambda(p){\ge}0~\nu{-}a.e.$ Let $T$ be a countable subset of $T_f$ which is dense in $T_f$ with respect to the finest locally convex topology on $\dC[x]$ (for instance, one may take polynomials $g$ in the definition of $T_f$ with rational coefficients). For each $t\in T$ there  is $\nu$-null subset $N_t\supseteq N$ of $\sigma (H)$ such that $L_\lambda(t)\ge 0$ for $\lambda\in\sigma(H)/N_t$. Setting $N^\prime {:=}\cup_tN_t$, we have $L_\lambda(T)\ge 0$ and hence $L_\lambda(T_f)\ge 0$ for all $\lambda\in\sigma(H)/N^\prime.$ The relation $L_\lambda((h_j{-}\lambda_j)p){=}0$ for $\lambda\in\sigma(H)/N$ is obvious by the definition of $L_\lambda$.\hfill$\Box$

\sn
In the proof of Theorem 1 we use Haviland's theorem [H] which is stated as

\mn
{\bf Lemma 7.} {\it Let $M$ be a closed subset of $\dR^d$ and let $F$ be a linear functional on $\dC[x]$ such that $F(p){\ge}0$ whenever $p\in\dR[x]$ and $p\ge 0$ on $M$. Then there exists a positive Borel measure $\mu$ on $M$ such that $F(p){=}\int_M p d\mu $ for all $p\in\dC[x]$.}

\mn
{\bf Proof of Theorem 1:} Let $L$ be a linear functional on $\dC[x]$ such that $L(T_f){\ge}0$. Suppose that for all $\lambda\in\Lambda$ the set $K_f\cap C_\lambda$ has property $(SMP)$. Fix a polynomial $p\in\dR[x]$ such that $p\ge 0$ on $K_f$. Let $\lambda\in\sigma(H)/N^\prime$. Then, $L_\lambda(T_{f(\lambda)} )\ge 0$ by Lemma 6(ii). Since the set $K_f\cap C_\lambda$ has property $(SMP)$, the functional $L_\lambda$ comes from a positive Borel measure on $K_f\cap C_\lambda$. Therefore, since $p\ge 0$ on $K_f$, we have $L_\lambda(p)\ge 0$. From Lemma 6(i) we obtain $L(p)\ge 0$. Thus, by Haviland's theorem (Lemma 7) applied to $M=K_f$, there is a positive Borel measure $\mu$ on $K_f$ such that $L(p)=\int p d\mu$ for $p\in\dC[x]$. This proves that $K_f$ has property $(SMP)$.

The proof for property $(MP)$ is given by a slight modification of the above reasoning. One takes a non-negative polynomial $p$ on $\dR^d$ and applies Haviland's theorem to $M=\dR^d$.\hfill$\Box$

\mn
{\it Remark.} The preceding proof shows a slightly stronger assertion. It suffices to assume the properties $(SMP)$ and $(MP)$ for $\lambda$ in $\sigma(H)$ rather than $\Lambda$. Also, it is sufficient that the functional $L_\lambda$ comes from a measure supported on $K_f$ rather than $K_f\cap C_\lambda$.\\

\mn
{\bf 4. Applications}

\sn
The following simple fact is used below.

\sn
{\bf Lemma 8.} {\it Let $K_f$ be a semi-algebraic set in $\dR^d$. If $K_f$ is a real line (resp. a subset of a real line), then $f$ has property $(SMP)$ (resp. $(MP)$).}

\sn
{\bf Proof.} Applying an affine transformation if necessary, we can assume that the line is $\dR\cdot x_1$. By Corollary 4 applied to $p{=}x_j, j{=}2,{\dots}, n$, we have $L(p(x_1,{\dots},x_d))=L(p(x_1,0,{\dots},0))$ for $p\in\dR[x]$. Hence the assertion follows from Hamburger's theorem.\hfill $\Box$

\sn
In order to apply Theorem 1 one has to look for sequences $f$ and bounded polynomials $h_1,{\dots},h_n$ on the semi-algebraic set $K_f$ such that for each $\lambda\in\Lambda$ the sequence $f(\lambda)$ of the ``fiber set"´ $K_f\cap C_\lambda$ has property $(SMP)$ or $(MP)$. In particular, $f$ has property $(SMP)$ if for each $\lambda\in\Lambda$ the set 
$K_{f(\lambda)}=K_f\cap C_\lambda$ is either compact (by [S1]) or a real line (by Lemma 8). If $K_f$ is a half-line, then $f$  does not have property $(SMP)$ in general (see Example 5 below). Clearly, the sequence $f$ has property $(MP)$ if $K_\lambda\cap C_\lambda$ is compact or a subset of a real line for all $\lambda\in \Lambda$.

Combined with the results of [S], Theorem 1 yields a large class of sequences having properties $(SMP)$ resp. $(MP)$.

\mn
{\bf Theorem 9.} {\it Let $h_1,{\dots},h_n,n{\in}\dN$, be polynomials of $\dR[x_1,{\dots},x_d]$ and let $m_1,{\dots},m_n, M_1,{\dots}$, $M_n{\in}\dR$ be such that $m_1{\le}M_1,{\dots},m_n{\le} M_n$. Let $f=\break
(f_1,{\dots},f_{2n})$ be the sequence given 
by $f_ {2j{-}1}{=}h_j{-}m_j, f_{2j}{=}M_j{-}h_j$, $j{=}1,{\cdots},n$. Suppose that for each $\lambda\in \Lambda\equiv[m_1,M_1]\times{\cdots}\times [m_n,M_n]$ the set $C_\lambda{:=}\{x\in\dR^d:h_j(x){=}\lambda_j,j{=}1,{\dots},n\}$ is empty or an irreducible smooth affine real curve which is rational or has at least one non-real point at infinity.

Then the defining sequence $f$ of the semi-algebraic set $K_f= \{x \in \dR^d : m_j \le h_j(x)$
$\le M_j, j = 1,{\dots} ,n\}$ has property $(SMP)$.}

\mn
{\bf Proof.} From Corollary 4 it follows that each linear functional $L$ on $\dC[x]$ such that $L(T_{f(\lambda)})\ge 0$ annihilates the vanishing ideal of the real algebraic variety $K_{f(\lambda)}=C_\lambda, \lambda\in \Lambda$. Therefore, by Theorems 3.11 and 3.12 in [PS] (which are essentially based on the results in [S]), the assumptions on $C_\lambda$ imply that the sequence $f(\lambda)$ has property $(SMP)$ for any $\lambda\in \Lambda$. Hence $f$ has property $(SMP)$ by Theorem 1.\hfill $\Box$

Perhaps the simplest tube sets $K_f$ in Theorem 9 are cylinders with compact base sets.

\mn
{\bf Corollary 10.} {\it Let $K$ be a compact semi-algebraic set in $\dR^d$. Let $f{=}(f_1,{\dots},f_k)$ be a sequence of polynomials $f_j\in \dR[x_1{\cdots},x_{d+1}]$. If the semi-algebraic set $K_f$
in $R^{d+1}$ is equal to $K\times \dR$, then $f$ has property $(SMP)$. If
$K_f$ is contained in $K\times \dR$, then $f$ has property $(MP)$.}

\mn
{\bf Proof.} Let $K=K_h$, where $h=(h_1,{\dots}, h_n)$ and $h_1,{\dots}, h_n\in\dR[x_1,{\dots}, x_d]$. Since $K$ is compact, the polynomials $h_j$ are bounded on $K_f$. Since the fiber set $K_f\cap C_\lambda$ is a subset of a real line, the assertion follows at once from Lemma 8 and Theorem 1.\hfill $\Box$

For the natural choice of generators of $K_f=K_h\times \dR$ the first assertion of Corollary 10 was also proved in [KM]. The very special case where $K$ is a ball was obtained in [Mc].

\sn
We illustrate our results by some examples.

\mn
{\bf Examples.} 1.) Let $m,M,c\in\dR$ be such that $c\ge 0$ and $0{<}m{\le} M$. 
Let $f{=}((x_1{-}c)x_2{-}m,M{-}(x_1{-}c)x_2)$. Then $h{:=}(x_1{-}c)x_2$ is bounded on the set
$K_f{=}\{x\in\dR^2{:}m{\le} h(x){\le} M\}$. For $\lambda\in[m,M]$ the sequence $f(\lambda)$ of the  fiber set $K_{f(\lambda)}{=}K_f\cap C_\lambda{=}\{x{:}(x_1{-}c)x_2{=}\lambda\}$ has property $(SMP)$. The latter follows at once from Example 3.7 in [KM] or from Theorem 3.12 in [PS]. Hence, by Theorem 1, $f$ has property $(SMP)$.\\
2.) The sets in this example and the next are subsets of the strip $[0,1]\times\dR$. Let $f{=}(x_1x_2,1{-}x_1x_2,x_1,1{-}x_1)$. Then the algebra $\cB(K_f)$ consists of all polynomials $q(x_1,x_1x_2)$, where $q\in\dC[x_1,x_2]$. Let $h_1=x_1x_2$ and $h_2=x_1$. Theorem 1 can be applied to $h_1$ or to $h_2$ or to $h_1$ and $h_2$. In all three cases we conclude that the sequence $f$ has property $(SMP)$. In the third case the fiber set $K_f\cap C_\lambda$ is a point if $\lambda_1\ne 0$ and the $x_2$-axis if $\lambda_1=0$.\\
3.) Let $f{=}(x_1x_2{-}1,2{-}x_1x_2,x_1,1{-}x_1)$. Applying Theorem 1 to $h_1=x_1x_2,$\break $h_2=x_1$, all fiber sets are points , so $f$ has property $(SMP)$.\\
4.) Theorem 1 can be also used to solve the moment problem on curves. Let $f=(f_1,{\dots}, f_k)$ be a sequence of polynomials $f_j\in\dR[x_1,{\dots}, x_d]$ such that $K_f$ is an irreducible real algebraic curve in $\dR^d$. Suppose that there exists a non-zero $(a_1,{\dots}, a_d)\in\dR^d$ such that the linear polynomial $h(x):=a_1x_1+{\cdots}+a_d  x_d$ is bounded on $K_f$. Then each fiber set $K_f\cap \{ x\in\dR^d:h(x)=\lambda\}$ is empty or compact. Therefore, by Theorem 1, the sequence $f$ has property $(SMP)$. Examples of such curves in $\dR^2$ are $x^3_1+x^3_2=1$ with $h(x)=x_1+x_2$ and $x^2_2(1-x_1)=x^3_1$ with $h(x)=x_1.$\\
5.) This example is mainly taken from [BM]. Let $d=1, f(x)=x^3$, so $K_f=[0,+\infty)$. Since each linear functional $L$ on $\dC[x]$ such that $L(T_f)\ge 0$ is non-negative on squares, by Hamburger's theorem it can be given by a positive measure on $\dR$. That is, $f$ has property $(MP)$. We show that $K_f$ does not obey property $(SMP)$.\\
Let $\mu$ be an $N$-extremal measure on $\dR$ such that $\mu$ represents an indeterminate Stieltjes moment sequence, ${\rm supp}~\mu ~ \cap (-\delta,\delta)=\emptyset$ for some $\delta>0$ and ${\rm supp}~\mu ~ \cap (-\infty,0)\ne\emptyset$. (The existence of such a measure follows easily from the existence of indetermiante Stieltjes moment sequences (cf. [ST], p. 22) combined with basic properties of Nevanlinna's extremal measures as collected in Theorem 2.13 in [ST].)
 Define $L(p)=\int px^{-2} d\mu, ~p\in\dC[x]$. Since $\mu$ gives a Stieltjes moment sequence, $L(T_f)\ge 0$. Since $\mu$ is $N$-extremal, $\dC[x)$ is dense in 
$L^2(\mu)$ ([A], Theorem 2.3.3). Therefore, $(x\pm i)\dC[x]$ is dense in $L^2(x^{-2}\mu)$ and hence the measure $x^{-2}d\mu$ is determinate.  Since this measure is not supported on $K_f=[0,+\infty)$, $L$ cannot be represented by a positive measure on $\dR$, that is, $f$ does not have property $(SMP)$. Note that there is a polynomial $p_0\in\dC[x]$ such that $L(x p_0\overline{p_0})<0$. (Otherwise, by Stieltjes theorem $L$ could be given by a positive measure on $(0,+\infty)$.)

Let $\tilde{f} (x)=x$. Then $K_{\tilde{f}}= K_f$ and ${\tilde{f}}$ has property $(SMP)$. This  shows that, in contrast to the compact case [S1], property $(SMP)$ depends in general on the defining polynomials $f$ rather than the set $K_f$.\hfill$\Box$

\mn
It is well-known that there exists a linear functional $L$ on $\dC[x_1,x_2]$ such that $L(p\bar{p})\ge 0$ for $p\in\dC[x_1,x_2]$ which is not a moment functional (that is, $L$ cannot be represented by a positive Borel measure on $\dR^2$). An explicit example can be found in [F]; it is reproduced in [S2], p. 62. We close this paper by constructing a much simpler explicit example which is even tractable for computations. It is based on the fact that no defining sequence of the curve $x^3_1=x^2_2$ in $\dR^2$ has property $(SMP)$. 

\mn
{\bf Example} 6.)~ Suppose that $L$ is a linear functional on $\dC[x]$ such that\\
(i) $L(p\bar{p} + x^3 q\bar{q})\ge 0$ for $p,q\in\dC[x]$.\\
(ii) There is a polynomial $p_0\in\dC[x]$ such that $L(xp_0\overline{p_0})<0$.

An {\it explicit} example of such a functional $L$ was derived in Example 5. By some slight modifications of the construction therein one can have that $L(x)<0$. The {\it existence} of such a functional follows also from the Hahn-Banach separation theorem of convex sets:  Since the preorder $T_f=\{p\bar{p}+x^3q\bar{q};p,q\in\dC[x]\}$ is closed in the finest locally convex toplogy on $\dC[x]$ and $x\notin T_f$ as easily shown, there is a linear functional $L$ on $\dC[x]$ such that $L(T_f)\ge 0$ and $L(x)<0$.

Define linear functionals $L_1$ on $\dC[x]$ and $L_2$ on $\dC[x_1, x_2]$ by $L_1(x^{2n})=L(x^n), L_1(x^{2n+1})=0$ for $n\in\dN_0$ and $L_2(p(x_1,x_2))=L_1(p(x^2,x^3))$ for $p\in\dC[x_1,x_2]$. Using the fact that $L(T_f)\ge 0$ one verifies that $L_2(p\bar{p})\ge 0$ for all $p\in\dC[x_1,x_2].$ But $L_2$ is not a moment functional. Assume to the contrary that $L_2$ is given by a positive Borel measure $\mu$ on $\dR^2$. By the definition of $L_2$ the support of $\mu$ is contained in the curve $x^3_1=x^2_2$. Since $x_1\ge 0$ on this curve and $L_2(x_1p_0(x_1)\overline{p_0(x_1)} ))=L_1(x^2p_0(x^2)\overline{p_0(x^2)} )=L(xp_0(x)\overline{p_0(x)})<0$, we get a contradiction. \hfill$\Box$

\mn
{\it Remarks.} The results of this paper have shown how the existence of ``sufficiently many" {\it bounded} polynomials on a closed semi-algebraic set $K_f$ can be used for the study of the moment problem on $K_f$. One might try to use such an assumption also for other problems such as the strict Positivstellensatz. Consider  the following property of a sequence $f$:\\
(SPS): {\it For each $p\in\dR[x]$ such that $p{\ge}0$ on $K_f$ there exists $q{\in} T_f$  such that \\
${}$\hspace{1.3cm}$p{+}\varepsilon q\in T_f$ for all $\varepsilon>0$.}\\
Clearly, $(SPS)$ implies $(SMP)$. As shown in [KM], $(SPS)$ holds for the natural choice of generators of a cylinder $K_h\times \dR$ with compact base set $K_h$. We conclude this paper by stating the following question: Does the assertion of Theorem 1 remain valid if property $(SMP)$ is replaced by $(SPS)$?

\mn

\end{document}